\input amstex
\documentstyle{amsppt}
\magnification=\magstep1                        
\hsize6.5truein\vsize8.9truein                  
\NoRunningHeads
\loadeusm

\magnification=\magstep1                        
\hsize6.5truein\vsize8.9truein                  
\NoRunningHeads
\loadeusm

\topmatter

\title
The Mahler measure of the Rudin-Shapiro polynomials
\endtitle

\rightheadtext{the Mahler measure of the Rudin-Shapiro polynomials}

\author
Tam\'as Erd\'elyi
\endauthor

\address Department of Mathematics, Texas A\&M University,
College Station, Texas 77843 \endaddress

\email terdelyi\@math.tamu.edu
\endemail

\endthanks


\abstract
Littlewood polynomials are polynomials with each of their coefficients in
$\{-1,1\}$. A sequence of Littlewood polynomials that satisfies a remarkable 
flatness property on the unit circle of the complex plane is given by the
Rudin-Shapiro polynomials. It is shown in this paper that the Mahler
measure and the maximum modulus of the Rudin-Shapiro polynomials on the unit
circle of the complex plane have the same size. It is also shown  that the 
Mahler measure and the maximum norm of the Rudin-Shapiro polynomials have the 
same size even on not too small subarcs of the unit circle of the complex plane. 
Not even nontrivial lower bounds for the Mahler measure of the Rudin Shapiro 
polynomials have been known before.
\endabstract

\endtopmatter

\document

\head 1. Introduction \endhead

Let $\alpha < \beta$ be real numbers. The Mahler measure $M_{0}(Q,[\alpha,\beta])$ is
defined for bounded measurable functions $Q$ defined on $[\alpha,\beta]$ as
$$M_{0}(Q,[\alpha,\beta]) := \exp\left(\frac{1}{\beta - \alpha} \int_{\alpha}^{\beta}
{\log|Q(e^{it})|\,dt} \right)\,.$$
It is well known, see [17], for instance, that
$$M_{0}(Q,[\alpha,\beta]) = \lim_{q \rightarrow 0+}{M_{q}(Q,[\alpha,\beta])}\,,$$
where
$$M_{q}(Q,[\alpha,\beta]) := \left( \frac{1}{\beta-\alpha} \int_{\alpha}^{\beta}
{\left| Q(e^{it}) \right|^q\,dt} \right)^{1/q}\,, \qquad q > 0\,.$$
It is a simple consequence of the Jensen formula that
$$M_0(Q) := M_0(Q,[0,2\pi]) = |c| \prod_{k=1}^n{\max\{1,|z_k|\}}$$
for every polynomial of the form
$$Q(z) = c\prod_{k=1}^n{(z-z_k)}\,, \qquad c,z_k \in {\Bbb C}\,.$$
See [3, p. 271] or [2, p. 3], for instance. 
Let $D := \{z \in {\Bbb C}\}: |z| < 1\}$ denote the open unit disk of the complex plane.
Let $\partial D :=  \{z \in {\Bbb C}\}: |z| = 1\}$ denote the unit circle of the complex plane.

Finding polynomials with suitably restricted coefficients and
maximal Mahler measure has interested many authors. The classes
$${\Cal L}_n :=\left\{ p: p(z) =\sum_{k=0}^{n}{a_kz^k}\,, 
\quad a_k \in \{-1,1\} \right\}$$
of Littlewood polynomials and the classes
$${\Cal K}_n := \left\{ p:p(z) = \sum_{k=0}^{n}{a_kz^k}\,, \quad a_k \in {\Bbb C}, 
\enskip |a_k| =1 \right\}$$
of unimodular polynomials are two of the most important classes considered.
Observe that ${\Cal L}_n \subset {\Cal K}_n$ and 
$$M_0(Q) = M_0(Q,[0,2\pi]) \leq M_2(Q,[0,2\pi]) = \sqrt{n+1}$$ 
for every $Q \in {\Cal K}_n$. 
Beller and Newman [1] constructed unimodular polynomials $Q_n \in {\Cal K}_n$  whose
Mahler measure $M_0(Q,[0,2\pi])$ is at least $\sqrt{n}-c/\log n$.

Section 4 of [2] is devoted to the study of Rudin-Shapiro polynomials. 
Littlewood asked if there were polynomials $p_{n_k} \in {\Cal L}_{n_k}$ satisfying  
$$c_1 \sqrt{n_k+1}  \leq |p_{n_k}(z)| \leq c_2 \sqrt{n_k+1}\,, \qquad z \in \partial D\,,$$
with some absolute constants $c_1 > 0$ and $c_2 > 0$, see [2, p. 27] for a reference 
to this problem of Littlewood.
To satisfy just the lower bound, by itself, seems very hard, and no such sequence $(p_{n_k})$  
of Littlewood polynomials $p_{n_k} \in {\Cal L}_{n_k}$ is known. A sequence of Littlewood polynomials 
that satisfies just the upper bound is given by the Rudin-Shapiro polynomials. The Rudin-Shapiro 
polynomials appear in  Harold Shapiro's 1951 thesis [20] at MIT and are sometimes called just 
Shapiro polynomials. They also arise independently in Golay's paper [16]. They are 
remarkably simple to construct and are a rich source of counterexamples to possible 
conjectures.

The Rudin-Shapiro polynomials are defined recursively as follows:
$$P_0(z) :=1\,, \qquad Q_0(z) := 1\,,$$
and
$$\split P_{n+1}(z) & := P_n(z) + z^{2^n}Q_n(z)\,, \cr
Q_{n+1}(z) & := P_n(z) - z^{2^n}Q_n(z)\,, \cr \endsplit$$
for $n=0,1,2,\ldots\,.$ Note that both $P_n$ and $Q_n$ are polynomials of degree $N-1$ with $N := 2^n$ 
having each of their coefficients in $\{-1,1\}$. 
It is well known and easy to check by using the parallelogram law that
$$|P_{n+1}(z)|^2 + |Q_{n+1}(z)|^2 = 2(|P_n(z)|^2 + |Q_n(z)|^2)\,, \qquad z \in \partial D\,.$$
Hence
$$|P_n(z)|^2 + |Q_n(z)|^2 = 2^{n+1} = 2N\,, \qquad z \in \partial D\,. \tag 1.1$$
It is also well known (see Section 4 of [2], for instance), that 
$$|Q_n(z)| = |P_n(-z)|\,, \qquad z \in \partial D\,. \tag 1.2$$
Peter Borwein's book [2] presents a few more basic results on the Rudin-Shapiro 
polynomials. Cyclotomic properties of the Rudin-Shapiro polynomials are discussed in 
[8]. Obviously $M_2(P_n,[0,2\pi]) = 2^{n/2}$ by the Parseval formula. In 1968 Littlewood 
[18] evaluated $M_4(P_n,[0,2\pi])$ and found that $M_4(P_n,[0,2\pi]) \sim (4^{n+1}/3)^{1/4}$. 
Rudin-Shapiro like polynomials in $L_4$ on the unit circle of the complex 
plane are studied in [6]. In 1980 Saffari conjectured that 
$$M_q(P_n,[0,2\pi]) \sim \frac{2^{(n+1)/2}}{(q/2+1)^{1/q}}$$  
for all even integers $q > 0$, perhaps for all real $q > 0$.
This conjecture was proved for all even values of $q \leq 52$ by Doche [11] 
and Doche and Habsieger [12].

Despite the simplicity of their definition not much is known
about the Rudin-Shapiro polynomials. It is shown in this paper that the Mahler
measure and the maximum modulus of the Rudin-Shapiro polynomials on the unit
circle of the complex plane have the same size.  A consequence of this result is also proved. 
It is also shown in this paper that the Mahler measure and the maximum norm of the 
Rudin-Shapiro polynomials have the same size even on not too small subarcs of the 
unit circle of the complex plane. Not even nontrivial lower bounds for the Mahler 
measure of the Rudin-Shapiro polynomials has been known before.

P. Borwein and Lockhart [7] investigated the asymptotic behavior of the mean
value of normalized $L_p$ norms of Littlewood polynomials for arbitrary $p > 0$.
They proved that
$$\lim_{n \rightarrow \infty} {\frac{1}{2^{n+1}} \, \sum_{f \in {\Cal L}_n}
{\frac{(M_q(f,[0,2\pi]))^q}{n^{q/2}}}}= \Gamma \left( 1+ \frac q2 \right)\,.$$
An analogue of this result does not seem to be known for $q=0$ (the Mahler measure).
However, the recent paper [9] paper establishes beautiful results on the average Mahler 
measure and $L_q({\partial D})$ norms of unimodular polynomials (polynomials with complex coefficients of 
modulus $1$) of degree $n$.

\medskip

\head 2 Main Theorems \endhead 
Our first theorem states that the Mahler measure and the maximum norm of the Rudin-Shapiro polynomials
on the unit circle of the complex plane have the same size.

\proclaim{Theorem 2.1}
Let $P_n$ and $Q_n$ be the $n$-th Rudin-Shapiro polynomials defined in Section 1. 
There is an absolute constant $c_1 > 0$ such that
$$M_0(P_n,[0,2\pi]) = M_0(Q_n,[0,2\pi]) \geq c_1\sqrt{N}\,,$$ 
where
$$N := 2^n = \text{\rm deg}(P_n)+1  = \text{\rm deg}(Q_n)+1\,.$$
\endproclaim

By following the line of our proof, it is easy to verify that $c_1 = e^{-227}$ is an appropriate choice 
in Theorem 2.1. To formulate our next theorem we define
$$\widetilde{P}_n := 2^{-(n+1)/2}P_n \qquad \text{and} \qquad \widetilde{Q}_n := 2^{-(n+1)/2}Q_n\,. \tag 2.1$$ 
By using the above normalization, (1.1) can be rewritten as 
$$|\widetilde{P}_n(z)|^2 + |\widetilde{Q}_n(z)|^2 = 1\,, \qquad z \in \partial D\,. \tag 2.2$$ 
Let
$$I_q(\widetilde{P}_n) := (M_q(\widetilde{P}_n,[0,2\pi]))^q := 
\frac{1}{2\pi} \int_0^{2\pi} {|\widetilde{P}_n(e^{i\tau})|^q \, d\tau}\,, \qquad q > 0\,.$$
The following result on the moments of the Rudin-Shapiro polynomials is a simple consequence 
of Theorem 2.1.

\proclaim{Theorem 2.2}
There is a constant $L < \infty$ independent of $n$ such that
$$\sum_{k=1}^\infty{\frac{I_k(\widetilde{P}_n)}{k}} < L\,, \qquad n=0,1,\ldots\,.$$ 
\endproclaim

Our final result states that the Mahler measure and the maximum norm of the Rudin-Shapiro 
polynomials have the same size even on not too small subarcs of the unit circle of the complex plane.

\proclaim{Theorem 2.3} There is an absolute constant $c_2 > 0$ such that
$$M_0(P_n,[\alpha,\beta]) \geq c_2 \sqrt{N}\,, \qquad N := 2^n = \text{\rm deg}(P_n)+1\,,$$
for all $n \in {\Bbb N}$ and for all $\alpha, \beta \in {\Bbb R}$ such that
$$\frac{32\pi}{N} \leq \frac{(\log N)^{3/2}}{N^{1/2}} \leq \beta - \alpha \leq 2\pi\,.$$
The same upper bound holds for $M_0(Q_n,[\alpha,\beta])$.
\endproclaim

It looks plausible that Theorem 2.3 holds whenever $32\pi/N \leq \beta - \alpha \leq 2\pi\,,$
but we do not seem to be able to handle the case
$32\pi/N \leq \beta - \alpha \leq (\log N)^{3/2}N^{-1/2}$ in this paper.

\medskip

\head 3. Lemmas \endhead 

A key to the proof of Theorem 2.1 is the following observation which is a straightforward 
consequence of the definition of the Rudin-Shapiro polynomials $P_n$ and $Q_n$.

\proclaim{Lemma 3.1} Let $n \geq 2$ be an integer, $N := 2^n$, and  and let
$$z_j := e^{it_j}\,, \quad t_j := \frac{2\pi j}{N}\,, \quad j \in {\Bbb Z}\,.$$ 
We have
$$\split P_n(z_j) & = 2P_{n-2}(z_j)\,, \qquad j=2u\,, \enskip u \in {\Bbb Z}\,, \cr
P_n(z_j) & =(-1)^{(j-1)/2} 2i\,Q_{n-2}(z_j)\,, \qquad j=2u+1\,, \enskip u \in {\Bbb Z}\,, \cr \endsplit$$ 
where $i$ is the imaginary unit.
\endproclaim

Another key to the proof of Theorem 2.1 is Theorem 1.3 from [15].
Let ${\Cal P}_N$ be the set of all polynomials of degree at most $N$ with real coefficients.

\proclaim{Lemma 3.2}
Assume that $n,m \geq 1$,
$$0 <\tau_1 \leq \tau_2 \leq \cdots \leq \tau_m \leq 2\pi\,, \quad \tau_0 := \tau_m - 2\pi\,, \quad \tau_{m+1} := \tau_1 + 2\pi\,.$$ 
Let
$$\delta := \max \{\tau_1 - \tau_0, \tau_2 - \tau_1, \ldots, \tau_m - \tau_{m-1}\}\,.$$
For every $A > 0$ there is a $B > 0$ depending only on $A$ such that
$$\sum_{j=1}^m{\frac{\tau_{j+1} - \tau_{j-1}}{2} \, \log |P(e^{i\tau_j})|}
\leq \int_0^{2\pi}{\log |P(e^{i\tau})| \, d\tau} + B$$
for all $P \in {\Cal P}_N$ and $\delta \leq AN^{-1}$. Moreover, the choice $B = 9A^2$ is appropriate. 
\endproclaim 

Our next lemma can be proved by a routine zero counting argument.
Let ${\Cal T}_k$ be the set of all real trigonometric polynomials of degree at most $k$.

\proclaim{Lemma 3.3} 
For $k \in {\Bbb N}$, $ M > 0$, and $\alpha \in {\Bbb R}$, let $T \in {\Cal T}_k$ be defined by
$$T(t) = \frac M2\,(1 - \cos(k(t - t_0))) = M\,\sin^2\left(\frac{k(t-t_0)}{2}\right)\,. \tag 3.1$$
Let $a \in {\Bbb R}$ be fixed. Assume that $S \in {\Cal T}_k$ satisfies $S(a) = T(a) > 0$ and
$0 \leq S(t) \leq M$ holds for all $t \in {\Bbb R}$.
Then

\noindent (i) $S(t) > T(t)$ holds for all $t \in (y,a)$ if $T$ is increasing on $(y,a)$.

\noindent (ii) $S(t) > T(t)$ holds for all $t \in (a,y)$ if $T$ is decreasing on $(a,y)$.
\endproclaim

\demo{Proof of Lemma 3.3}
If the lemma were false then $S-T \in {\Cal T}_k$ would have at least $2k+1$ zeros 
in a period, by counting multiplicities.
\qed \enddemo

A straightforward consequence of Lemma 3.3 is the following. For the sake of brevity let 
$$\gamma := \sin^2(\pi/8) = \frac 12 (1-\cos(\pi/4))\,.$$

\proclaim{Lemma 3.4}
Assume that $S \in {\Cal T}_k$ satisfies $0 \leq S(t) \leq M$ for all $t \in {\Bbb R}$.
Let $a \in {\Bbb R}$ and assume that  $S(a) \geq (1 - \gamma)M$. Then
$$S(t) \geq \gamma M\,, \quad t \in \left[a - \delta, a + \delta \right]\,, \quad \delta := \frac{\pi}{2k}\,.$$
\endproclaim

\demo{Proof of Lemma 3.4} 
Let $T \in {\Cal T}_k$ be defined by (3.1). Pick a $t_0 \in {\Bbb R}$ such that $S(a) = T(a)$ and $T^\prime(a) \geq 0$.
Now observe that $T$ is increasing on $[a-\delta,a+\delta]$ and $T(a-\delta) \geq \gamma M$, 
and Lemma 3.3 (i) gives the lower bound of the lemma for all $t \in [a-\delta,a]$. 
Now pick a $t_0 \in {\Bbb R}$ such that $S(a) = T(a)$ and $T^\prime(a) \leq 0$.
Now observe that $T$ is decreasing on $[a,a+\delta]$ and $T(a+\delta) \geq \gamma M$, 
and hence Lemma 3.3 (ii) gives the lower bound of the lemma for all $t \in [a,a+\delta]$.
\qed \enddemo

Combining Lemmas 3.1 and 3.4 we easily obtain the following.

\proclaim{Lemma 3.5}
Let $P_n$ and $Q_n$ be the $n$-th Rudin-Shapiro polynomials. Let $N:=2^n$ and $\gamma := \sin^2(\pi/8)$.
Let
$$z_j := e^{it_j}\,, \quad t_j := \frac{2\pi j}{N}\,, \qquad j \in {\Bbb Z}\,.$$
We have
$$\max \{|P_n(z_j)|^2,|P_n(z_{j+r})|^2\} \geq \gamma 2^{n+1} = 2\gamma N\,, \quad r \in \{-1,1\}\,,$$
for every $j=2u$, $u \in {\Bbb Z}$.
\endproclaim

Lemma 3.5 tells us that the modulus of the Rudin-Shapiro polynomials $P_n$ is certainly larger than $\sqrt{2\gamma N}$ 
at least at one of any two consecutive $N$-th root of unity, where $N := 2^n$. This is a crucial observation of this paper, 
and despite its simplicity it does not seem to have been observed before in the literature or elsewhere. 
Moreover, note that while our Theorems 2.1 and 2.3 are proved with rather small multiplicative positive absolute 
constants, Lemma 3.5 is stated with a quite decent explicit constant. 

\demo{Proof of Lemma 3.5}
Let $k := 2^{n-2}$, $j=2u$, $u \in {\Bbb Z}$. We introduce the trigonometric polynomial $S \in {\Cal T}_k$ by 
$$S(t) := |Q_{n-2}(e^{it})|^2\,, \qquad t \in {\Bbb R}\,.\tag 3.2$$
Recall that 
$$S(t) = |Q_{n-2}(e^{it})|^2 \leq 2^{n-1} \tag 3.2$$
Assume that 
$$|P_n(z_j)|^2 < \gamma 2^{n+1}\,.$$ 
Then (1.1) implies that
$$|Q_n(z_j)|^2  > (1-\gamma)2^{n+1}\,. \tag 3.3$$
By Lemma 3.1 we have 
$$|P_{n-2}(z_j)|^2 = \frac 14 \, |P_n(z_j)|^2\,,$$ 
and hence (1.1) implies that
$$|Q_{n-2}(z_j)|^2 = \frac 14 \, |Q_n(z_j)|^2\,.$$ 
Combining this with (3.3), we obtain 
$$|Q_{n-2}(z_j)|^2 = \frac 14 \, |Q_n(z_j)|^2 > (1-\gamma) 2^{n-1}\,.$$ 
Hence, using Lemma 3.4 with $S \in {\Cal T}_k$ and  $M := 2^{n-1}$ (recall that (3.2) and (1.1) imply that 
$0 \leq S(t) \leq 2^{n-1}$ for all $t \in {\Bbb R}$), we can deduce that 
$$|Q_{n-2}(z_{j+r})|^2 \geq \gamma 2^{n-1}\,, \qquad r \in \{-1,1\}\,.$$
Finally we use Lemma 3.1 again to conclude that
$$|P_n(z_{j+r})|^2 = 4|Q_{n-2}(z_{j+r})|^2 \geq \gamma 2^{n+1}\,, \qquad r \in \{-1,1\}\,.$$ 
\qed \enddemo

To prove Theorem 2.3 we need Theorem 2.1 from [13]. We state it as our next lemma 
by using a slightly modified notation.

\proclaim{Lemma 3.6} Let $\omega_1 < \omega_2 \leq \omega_1 + 2\pi\,,$
$$\omega_1 \leq \theta_1 < \theta_2 < \cdots < \theta_{\mu} \leq \omega_2\,,$$
$$\theta_0 := \omega_1 - (\theta_1 - \omega_1)\,, \qquad \theta_{\mu+1} := \omega_2 + (\omega_2 - \theta_{\mu})\,,$$
$$\delta := \max\{\theta_1 - \theta_0, \theta_2 - \theta_1,\ldots, \theta_{\mu+1} - \theta_{\mu}\} 
\leq \frac 12 \, \sin \frac{\omega_2 - \omega_1}{2}\,.$$
There is an absolute constant $c_3 > 0$ such that
$$\sum_{j=0}^{\mu}{\frac{\theta_{j+1}-\theta_{j-1}}{2} \log|P(e^{i\theta_j})|} 
\leq \int_{\omega_1}^{\omega_2}{\log|P(e^{i\theta})|\,d\theta} + c_3E(N,\delta,\omega_1,\omega_2)$$
for every polynomial $P$ of the form
$$P(z) = \sum_{j=0}^N{b_j z^j}\,, \qquad b_j \in {\Bbb C}\,, \enskip b_0b_N \neq 0\,,$$
where
$$E(N,\delta,\omega_1,\omega_2) := (\omega_2-\omega_1)N\delta + N\delta^2\log(1/\delta) + 
\sqrt{N \log R} \left(\delta \log(1/\delta) + \frac{\delta^2}{\omega_2-\omega_1} \right)$$
and $R := |b_0b_N|^{-1/2} \|P\|_{\partial D}\,.$
\endproclaim

Observe that $R$ appearing in the above theorem can be easily estimated by
$$R \leq |b_0b_N|^{-1/2}(|b_0| + |b_1| + \cdots + |b_N|)\,.$$

\medskip

\head 4. Proof of Theorems 2.1, 2.2, and 2.3 \endhead

\demo{Proof of Theorem 2.1}
Let, as before, $\gamma = \sin^2(\pi/8)$, $N := 2^n$, and 
$$z_j := e^{it_j}\,, \quad t_j := \frac{2\pi j}{N}\,, \quad j \in {\Bbb Z}\,.$$
By Lemma 3.5 we can choose 
$$\tau_m - 2\pi =: \tau_0 \leq 0 < \tau_1 < \tau_2 < \cdots < \tau_m \leq 2\pi$$ 
so that 
$$\{\tau_1,\tau_2,\ldots,\tau_m\} \subset \{t_1,t_2,\ldots,t_N\}\,, \tag 4.1$$
$$\tau_{j+1}-\tau_j \leq \frac{4\pi}{N}\,, \qquad j=0,1,\ldots,m-1\,, \tag 4.2$$
and  
$$|P_n(e^{i\tau_j})|^2 \geq \gamma 2^{n+1}\,, \qquad j=1,2,\ldots, m\,. \tag 4.3$$ 
Then the value 
$$\delta := \max \{\tau_1 - \tau_0, \tau_2 - \tau_1, \ldots, \tau_m - \tau_{m-1}$$
appearing in Lemma 3.2 satisfies $\delta \leq AN^{-1}$  with $A=4\pi$.
Let $B > 0$ be chosen for $A := 4\pi$ according to Lemma 3.2. Combining $P_n \in {\Cal P}_N$, (4.1), and Lemma 3.2, 
we conclude that
$$\split 2\pi \left( \frac 12 \log 2^{n+1} + \frac 12 \log \gamma \right) & \leq 
\sum_{j=1}^m{\frac{\tau_{j+1} - \tau_{j-1}}{2} \, \log |P_n(e^{i\tau_j})|} \cr
& \leq \int_0^{2\pi}{\log |P_n(e^{i\tau})| \, d\tau} + B\,, \cr \endsplit$$
and hence
$$M_0(P_n,[0,2\pi]) \geq \exp(-B/(2\pi))\sqrt{2\gamma}\,2^{n/2} = c\sqrt{N}$$
follows with the absolute constant $c := \exp(-B/(2\pi))\sqrt{2\gamma} > 0$. 
Combining this with (1.2), we obtain
$$M_0(Q_n,[0,2\pi]) \geq c2^{(n+1)/2}$$
with the same absolute constant $c > 0$.
\qed \enddemo

\demo{Proof of Theorem 2.2}
Recalling (2.2) and using the power series expansion of the function $f(z) : = \log(1-z)$ on $(-1,1)$, 
and the Monotone Convergence Theorem, we deduce that
$$\split \int_0^{2\pi}{\log|\widetilde{P}_n(e^{i\tau})|^2 \, d\tau} & = \int_0^{2\pi}{\log(1-|\widetilde{Q}_n(e^{i\tau})|^2) \, d\tau}  
= \int_0^{2\pi}{-\sum_{k=1}^\infty{\frac{|\widetilde{Q}_n(e^{i\tau})|^{2k}}{k} \, d\tau}} \cr
& = \int_0^{2\pi}{-\sum_{k=1}^\infty{\frac{|\widetilde{P}_n(e^{i\tau})|^{2k}}{k} \, d\tau}} 
= -\sum_{k=1}^\infty{\int_0^{2\pi}{\frac{|\widetilde{P}_n(e^{i\tau})|^{2k}}{k} \, d\tau}} \cr 
& = -2\pi \sum_{k=1}^\infty{\frac{I_{2k}(\widetilde{P}_n)}{k}}\,.  \cr \endsplit$$
Combining this with Theorem 2.1 gives that there is an $L < \infty$ independent of $n$ such that 
$$\sum_{k=1}^\infty{\frac{I_{2k}(\widetilde{P}_n)}{k}} < L\,.$$
As $I_k(\widetilde{P}_n)$ is a decreasing function of $k \in {\Bbb N}$, the theorem follows. 
\qed \enddemo

\demo{Proof of Theorem 2.3}
The theorem follows from Lemmas 3.5 and 3.6 in a straightforward fashion.
Note that  
$$(M_0(f,[\alpha,\beta]))^{\beta-\alpha} = (M_0(f,[\alpha,\gamma]))^{\gamma-\alpha}(M_0(f,[\gamma,\beta]))^{\beta-\gamma}\,,$$ 
for all $\alpha < \gamma < \beta \leq \alpha + 2\pi$ and for all functions $f$ continuous on $[\alpha,\beta]$. Hence, to prove the 
theorem, without loss of generality we may assume that $\beta-\alpha \leq \pi$.  
Let, as before, $\gamma = \sin^2(\pi/8)$, $N := 2^n$, and
$$z_j := e^{it_j}\,, \quad t_j := \frac{2\pi j}{N}\,, \quad j \in {\Bbb Z}\,.$$
By Lemma 3.5 we can choose
$$\tau_m - 2\pi =: \tau_0 \leq 0 < \tau_1 < \tau_2 < \cdots < \tau_m \leq 2\pi$$
so that (4.1), (4.2), and (4.3) hold. Let
$$\{\theta_1 < \theta_2 < \cdots < \theta_{\mu}\} := \{\tau_j \in [\alpha,\beta]: j=1,2,\ldots,m\}\,. \tag 4.4$$ 
The assumption on $N$  guarantees that the value of $\delta$ defined in Lemma 3.6 is at most $4\pi/N$ and 
$$\frac{4\pi}{N} \leq \frac{\beta-\alpha}{8} \leq \frac 12 \sin \frac{\beta - \alpha}{2}\,.$$ 
Observe also that $P_n \in {\Cal L}_{N-1}$, and hence when we apply Lemma 3.6 to $P_n$ we have $R \leq N$. 
By (4.1) we have 
$$|P_n(e^{i\theta_j})|^2 \geq \gamma 2^{n+1}\,, \qquad j=0,1,\ldots, \mu\,.$$
Applying Lemma 3.6 with $P := P_n$, $N := 2^n$, and $\{\theta_1 < \theta_2 < \cdots < \theta_{\mu}\}$ defined by (4.4) 
we obtain
$$\split (\beta - \alpha)\left( \frac 12 \log 2^{n+1} + \frac 12 \log \gamma \right) & \leq  
\sum_{j=0}^{\mu}{\frac{\theta_{j+1} - \theta_{j-1}}{2} \log|P_n(e^{i\theta_j})|} \cr 
& \leq \int_{\alpha}^{\beta}{\log|P_n(e^{i\theta})|\,d\theta} + c_3E(N,4\pi/N,\alpha,\beta)\,,\cr \endsplit$$
where the assumption 
$$\frac{(\log N)^{3/2}}{N^{1/2}} \leq \beta - \alpha \leq 2\pi$$ 
implies that 
$$\split E(N,4\pi/N,\alpha,\beta) \leq & \, c_4 \left( \frac{(\beta - \alpha)N}{N}
+ \frac{\log N}{N} + \sqrt{N\log N}\left(\frac{\log N}{N} + \frac{1}{N^2(\beta-\alpha)} \right) \right) \cr
\leq & \, c_5(\beta-\alpha) \cr \endsplit$$ 
with absolute constants $c_4 > 0$ and $c_5 > 0$. Hence  
$$M_0(P_n,[\alpha,\beta]) \geq \exp(-c_3c_5)\,\sqrt{2\gamma}\,2^{n/2} = c\sqrt{N}$$
with the absolute constant $c := \exp(-c_3c_5)\,\sqrt{2\gamma} > 0$. Combining this with (1.2), we obtain 
$$M_0(Q_n,[\alpha,\beta]) \geq c\sqrt{N}$$
with the same absolute constant $c > 0$.
\qed \enddemo

\medskip

\head 5. Some Recent results on Fekete Polynomials \endhead

This section could be viewed as part of the introduction. As the results on Fekete 
polynomials mentioned  in this section are somewhat analogous to our new results on 
the Rudin-Shapiro polynomials, we close this paper with this short section.  

For a prime number $p$ the $p$-th Fekete polynomial is defined as
$$f_p(z) := \sum_{k=1}^{p-1}{\left( \frac kp \right)z^k}\,,$$
where
$$\left( \frac kp \right) =
\cases
1, \quad \text{if \enskip} x^2 \equiv k \enskip (\text {mod\,}p) \enskip
\text{has a nonzero solution,}
\\
0, \quad \text{if \enskip} p \enskip \text{divides} \enskip k\,,
\\
-1, \quad \text{otherwise}
\endcases$$
is the usual Legendre symbol. Since $f_p$ has constant coefficient $0$, it is not
a Littlewood polynomial, but $g_p$ defined by $g_p(z) := f_p(z)/z$ is a Littlewood
polynomial of degree $p-2$, and has the same Mahler measure as $f_p$. Fekete polynomials are examined
in detail in [2], [4], [5], [10]. In [19] Montgomery proved that there are absolute constants
$c_1 > 0$ and $c_2 > 0$ such that
$$c_1 \sqrt{p} \log \log p \leq \max_{z \in \partial D}{|f_p(z)|} \leq c_2 \sqrt{p} \log p\,.$$
In [15] we gave a lower bound for the Mahler measure of the Fekete polynomials.
Namely we showed that for every $\varepsilon > 0$ there is a constant $c_{\varepsilon}$
such that
$$M_0(f_p,[0,2\pi]) \geq \left(\frac 12 - \varepsilon \right)\sqrt{p}$$
for all primes $p \geq c_{\varepsilon}$. One of the key lemmas in the proof of this
theorem is a remarkable property of the Fekete polynomials observed by Gauss. It states that
$$|f_p(z_p^j)| = p^{1/2}\,, \qquad j=1,2,\ldots, p-1\,,$$
and $f_p(1) = 0$, where $z_p := \exp(2\pi i/p)$ is the first $p$-th root of unity.
A simple proof of it is given in [2, pp. 37-38].
In [13] Theorem 1.2 is extended to subarcs of the unit circle. It is shown in [E-11]
that there is an absolute constant $c_1 > 0$ such that
$$M_0(f_p,[\alpha,\beta]) \geq c_1 p^{1/2}$$
for all prime numbers $p$ and for all $\alpha, \beta \in {\Bbb R}$ such that
$(\log p)^{3/2}p^{-1/2} \leq \beta - \alpha \leq 2\pi$.
In [14] the author gave an upper bound for the average value of $|f_p(z)|^q$ over any subarc $I$
of the unit circle, valid for all sufficiently large primes $p$ and arbitrary real exponents $q > 0$.
Namely, there is a constant $c_2(q,\varepsilon)$ depending only on $q > 0$  and $\varepsilon > 0$
such that
$$M_q(f_p,[\alpha,\beta]) \leq c_2(q,\varepsilon)p^{1/2}\,,$$
for all prime numbers $p$  and for all $\alpha, \beta \in {\Bbb R}$ such that
$\beta - \alpha \geq 2p^{-1/2+\varepsilon}$.

\bigskip

\subhead {6. Acknowledgment}\endsubhead
The author thanks Peter Borwein, Stephen Choi, Michael Mossinghoff, and Bahman Saffari for their careful reading 
of the paper and their comments on making the paper better.

\enddocument